\documentclass[matzei]{svjour}[]

\input{xypic}
\usepackage{amscd,amsmath,amssymb,amsfonts}
\usepackage[cmtip, all]{xy}
\usepackage{oldgerm}
\usepackage{hhline}

\newcommand{\liminv}[1]{{\displaystyle{\mathop{\rm 
lim}_{\buildrel\longleftarrow\over{#1}}}}\,}

\def\Q{\mathbb{Q}}
\def\Z{\mathbb{Z}}

\def\cal{\fam2}
\def\OO{{\cal{O}}}
\def\R{\mathbb{R}}
\def\C{\mathbb{C}}

\def\cD{\cal D}
\def\cF{\cal{F}}
\def\cH{\cal{H}}

\def\gr{\operatorname{gr}}

\def\zar{\operatorname{Zar}}
\def\zarloc{\operatorname{ZL}}
\def\an{\operatorname{an}}
\def\Ext{\operatorname{Ext}}
\def\Hom{\operatorname{Hom}}
\def\im{\operatorname{im }}
\def\Spec{\operatorname{Spec}}
\def\Pic{\operatorname{Pic}}
\def\Spec{\operatorname{Spec}}

\def\NS{\operatorname{NS}}
\def\Ext{\text{Ext}}
\def\MHS{\operatorname{MHS}}

\def\coker{\operatorname{coker}}
\def\ker{\operatorname{ker}}

\def\ns{\operatorname{ns}}
\def\L{\operatorname {L}}
\def\Xs{X_{\bullet}}
\renewcommand{\d}{\bullet}
\def\Hb{\mathbb H}
\def\H{\operatorname H}

\def\W{\operatorname{W}}
\def\Eb{\mathbb{E}}
\def\e{\operatorname{e}}
\def\J{\operatorname{J}}
\def\Lb{\mathbb {L}}
\def\L{\operatorname L}
\def\F{\operatorname F \hspace{-0.03cm}}
\def\E{\operatorname E}
\def\K{\operatorname K}
\def\Ker{\operatorname{Ker}}
\def\Cok{\operatorname{Coker}}
\def\Xsp{X'_{\d}}
\def\Rb{\mathbb R}
\def\di{\operatorname{d}}
\def\c{\operatorname{c}}
\def\T{\operatorname{T}}
\def\A{\operatorname{A}}
\def\Ib{\mathbb I}
\def\I{\operatorname{I}}

\def\R{\operatorname R}

\def\sn{\operatorname{sn}}

\begin{document}

\title{The N\'eron-Severi group of a proper seminormal complex 
variety}

\author{Luca~Barbieri-Viale\inst{1}\and 
Andreas Rosenschon\inst{2}\and 
V.~Srinivas\inst{3}}

\institute{Dipartimento di Matematica, Universit\`a di Padova, Padova, 
Italy,\\
\email{barbieri@math.unipd.it} 
\and Mathematisches Institut, Ludwigs-Maximilians-Universit\"at, M\"unchen, 
Germany, \\
\email{Andreas.Rosenschon@mathematik.uni-muenchen.de} 
\and Tata Institute of Fundamental Research, Mumbai, India, \\ 
\email{srinivas@math.tifr.res.in}}

\maketitle{}


\begin{abstract} 
We prove a Lefschetz $(1,1)$-Theorem for proper seminormal 
varieties over the complex numbers.
\end{abstract}

\section{Introduction}
\noindent
Let $X$ be a proper variety over $\C$. The N\'eron-Severi group 
$\NS(X)$ is the image of the map $\Pic(X)\rightarrow \H^2(X,\Z)$ 
which associates to a line bundle its cohomology class. If $X$ is 
smooth, the Lefschetz $(1,1)$-theorem \cite{Le} states 
$$
\NS(X)=\F^1\cap \H^2(X,\Z)=
\{x\in \H^2(X,\Z) \ | \ x_{\C}\in \F^1\H^2(X,\C)\}.
$$
In the smooth case one also knows the N\'eron-Severi group is 
isomorphic to the subgroup 
$$\H^2_{\zarloc}(X,\Z)\subseteq 
\H^2(X,\Z)$$ 
of Zariski locally trivial elements, which in turn is identified with
the Zariski cohomology group $\H^1_{\zar}(X, {\cH}^1_X)$, yielding a
formula
$$
\NS(X)= \H^2_{\zarloc}(X,\Z)=\H^1_{\zar}(X, {\cH}^1_X).
$$
\indent
If $X$ is singular, $\H^2(X,\Z)$ carries a mixed Hodge structure, 
and it makes sense to consider the question whether the above 
characterizations of the N\'eron-Severi group hold in this case as 
well; see Bloch's letter to Jannsen in \cite[Appendix A]{Jamotives}.
This is false. If $X$ is an irreducible  singular variety one always 
has inclusions
$$\displaylines{\NS(X)\subseteq
\H^2_{\zarloc}(X,\Z)=\H^1_{\zar}(X,{\cH}^1_X),\cr
\NS(X) \subseteq 
\F^1\cap \H^2(X,\Z)\cr}$$
and one knows of examples where both inclusions are strict 
\cite{BaSr1}, \cite{BaSr2}. Barbieri-Viale and Srinivas 
\cite[Question 4]{BaSr2}
raised the question whether the correct formulation of the Lefschetz 
theorem for a projective normal variety $X$ is 
\begin{equation}
\label{normallefschetz}
\NS(X)=\F^1\cap \H^1_{\zar}(X, {\cH}^1_X)=
\{x\in \H^1_{\zar}(X,{\cH}^1_X) \ | \ x_{\C}\in \F^1\H^2(X,\C)\}.
\end{equation}
This has been proved recently by Biswas-Srinivas \cite[Theorem 1.1]{BiSr}. 
They also give an example of a non-normal variety for which this Lefschetz 
theorem fails.
\\
\indent
Our main theorem states that the above description (\ref{normallefschetz}) 
also holds in the irreducible seminormal case. Recall that a variety is 
seminormal if all of its local rings are seminormal (see \cite{Sw} for 
the precise definition of a seminormal ring). Geometrically a variety $X$ 
is seminormal if any finite birational morphism $X'\rightarrow X$ that is 
bijective on closed points is an isomorphism.
For example, the nodal cubic 
curve is seminormal but the cuspidal cubic curve is not 
seminormal. For further discussion of the notion of seminormality from a 
scheme-theoretic viewpoint we refer to, for instance, \cite{GT}. 
\vspace{0.05cm}\\
\indent
We in fact prove the following more general result. 

\begin{theorem}\label{seminormallefschetz}
Let $X$ be a proper seminormal variety over $\C$. Then 
$$\NS(X)= {\F}^1 \cap \H^2_{\zarloc}(X, \Z).$$
\end{theorem}

Our methods differ considerably from the ones used in \cite{BiSr}. 
Instead of working with a resolution of singularities we use
simplicial techniques. Recall from \cite{DIII} that the mixed Hodge 
structure on $\H^2(X,\Z)$ is obtained from a smooth proper hypercovering 
$\pi_{\d}: \Xs \rightarrow X$ and the cohomological descent 
isomorphism $\H^*(X,\Z)\cong \Hb^*(\Xs,\Z_{\d})$. 
The basic idea behind our proof is to use such a hypercovering and 
cohomological descent to compare cohomology classes of line bundles on 
$X$ with cohomology classes of simplicial line bundles on $\Xs$. 
The assumption that $X$ is seminormal is equivalent to
$\pi_{\d*}{\OO}_{\Xs}\cong {\OO}_{X}$, and implies
$\Pic^0(X)\cong \Pic^0(\Xs)$, see Lemmas \ref{seminormal-simplicial} 
and \ref{exponentialseq}.
\\
\indent
Let $\Xs$ be a smooth proper simplicial scheme over $\C$. We study
the largest sub Hodge structure of $\Hb^2(\Xs,\Z_{\d})/\W_0$ (modulo 
torsion) of level $\leq 1$. By a standard argument \cite[\S 10]{DIII},
this Hodge structure $\Eb$ corresponds to a 1-motive $\Lb(\Eb)$, which is
given by an extension class map. It is a special case of the main result
of \cite{BRS} (see also \cite{Ra}) that the 1-motive $\Lb(\Eb)$ can be 
obtained from the Picard functors of the components $X_p$ of $\Xs$. 
More precisely, $\Eb$ is the Hodge realization of the isogeny class 
of the {\it Lefschetz} 1-motive $\Lb(\Xs)$ defined by the map
$$\ker\{\NS(X_0)\rightarrow \NS(X_1)\}\rightarrow \frac
{
\ker^0\{\Pic^0(X_1)\rightarrow \Pic^0(X_2)\}
}
{
\im\{\Pic^0(X_0)\rightarrow \Pic^0(X_1))\}
}$$
where $\ker^0$ denotes the connected component of the kernel of the
underlying homomorphism of algebraic groups.\\
\indent
We give an alternate characterization of the Hodge realization of 
$\Lb(\Xs)$. Let ${\cH}^1_{\d}(\Z)$ be the simplicial (Zariski) sheaf
whose $p$-th component is the sheaf associated to the presheaf
that sends a Zariski open subset $U\subseteq X_p$ to the (singular)
cohomology group $\H^1(U,\Z)$. It follows from a Leray spectral sequence 
argument that $\Hb^1_{\zar}(\Xs, {\cH}^1_{\d}(\Z))$ carries a mixed Hodge 
structure as a subgroup of $\Hb^2(\Xs,\Z_{\d})/\W_0$. We show 
$\Hb^1_{\zar}(\Xs, {\cH}^1_{\d}(\Q))$ is canonically isomorphic to 
$\Eb_{\Q}$; thus it is the Hodge realization of the isogeny class of 
$\Lb(\Xs)=\Lb(\Eb)$. 
\\
\indent
To prove our Lefschetz theorem we consider the simplicial Picard group
$$\Pic(\Xs)=\Hb^1_{\zar}(\Xs,\OO^*_{\Xs})=
\Hb^1_{\an}(\Xs, \OO^*_{\Xs}).$$
Let $\Pic^0(\Xs)$ be its connected component of the identity
and define the simplicial N\'eron-Severi group as the quotient 
$\NS(\Xs):=\Pic(\Xs)/\Pic^0(\Xs)$. The N\'eron-Severi group $\NS(\Xs)$ is the 
image of two cycle class maps, the analytic $\c_1: \Pic(\Xs)
\rightarrow \Hb^2(\Xs, \Z_{\d})$ given by the simplicial exponential sequence, 
and the map $\tilde{\c}_1: \Pic(\Xs) \rightarrow 
\Hb^1_{\zar}(\Xs, {\cH}^1_{\d}(\Z))$ defined by `pushing the exponential 
sequence to the Zariski site'. We compare the inclusions of the 
N\'eron-Severi group into the group of Zariski locally trivial classes 
on $X$ and $\Xs$ to obtain
$$\frac{\F^1\cap \H^2_{\zarloc}(X,\Z)}{\NS(X)} = \ker \{\delta:
\frac{\NS(\Xs)}{\NS(X)}\rightarrow 
\frac{\Hb^1_{\zar}(\Xs, {\cH}^1_{\d})}{ \im \H^2_{\zarloc}(X,\Z)}\}.
$$
If $X$ is seminormal we show $\ker \delta=0$. This last step involves
a non-trivial geometric argument (Lemma \ref{piclemma} below)
and does not seem to follow immediately from cohomological 
considerations. 
\\
\\
\noindent
{\bf Acknowledgments.} The first two authors thank Morihiko Saito for 
his comments about an earlier version of this paper. They also gratefully 
acknowledge financial support by the European Community TMR network {\it 
Algebraic $K$-theory, Linear Algebraic Groups and Related Structures}. The 
second author thanks TIFR for hospitality and excellent working conditions 
while this paper was written. 
\\
\\
\noindent
{\bf Notations.} If $X$ is a complex algebraic variety and $\cF$ is an 
analytic sheaf on $X$, we write $\H^p(X,\cF)$ for the corresponding 
cohomology groups; that is, unless specified otherwise, $\H^p$ denotes
cohomology in the analytic sense. Similarly, if $\Xs$ is a simplicial 
scheme and $\cF_{\d}$ is a simplicial sheaf, $\Hb^p(\Xs,\cF_{\d})$ denotes
the hypercohomology with respect to the analytic topology.


\section{Lefschetz 1-motive}
\noindent
In this section we define the Lefschetz 1-motive. Let $\Xs$ be a smooth 
proper simplicial scheme over $\C$. We are interested in the largest 
Hodge substructure $\Eb$ of 
$\im \{\Hb^2(\Xs,\Z)\rightarrow \Hb^2(\Xs,\Q_{\d})/{\W_0}\}$ 
of level $\leq 1$. This Hodge substructure is given by the extension
\begin{equation}\label{Ebdef}
0 \rightarrow \gr^{\W}_1 \Hb^2 \rightarrow \Eb \rightarrow
\Hb^{(1,1)}_{\Z}\rightarrow 0,
\end{equation}
where ${\Hb}^{(1,1)}:=({\gr}^{\W}_2 {\Hb}^2)^{(1,1)}$
\cite[2.3.7]{DII} and $\Eb$ is the pullback of $\Hb^{(1,1)}_{\Z}\subseteq 
{\gr}^{\W}_2\Hb^2$ along the projection map $\Hb^2/\W_0\rightarrow
{\gr}^{\W}_2\Hb^2$ \cite[2.3.7]{DIII}.
A standard argument \cite[\S 10]{DIII} shows that $\Eb$ is the 
Hodge realization of the 1-motive
\begin{equation}
\Lb(\Eb):= [\Hb^{(1,1)}_{\Z} \overset{\e}{\rightarrow}
\J^1({\gr}^{\W}_1\Hb^2)],
\end{equation}
where $\e$ is the extension class map. It is implicit in \cite{BRS} that the 
1-motive $\Lb(\Eb)$ is algebraically defined, see Proposition \ref{motiso} 
below for a precise statement. For the convenience of the reader we briefly 
sketch the argument.
\\
\indent
Let $\Xs$ be a proper simplicial scheme over an algebraically closed field $k$ 
of characteristic 0, and let $\Xsp$ be a smooth proper hypercovering of $\Xs$. 
On each component $X'_p$ we have an exact sequence
\begin{equation}
\label{componentpicext}
0\rightarrow \Pic^0(X'_p) \rightarrow \Pic(X'_p) \rightarrow
\NS(X'_p)\rightarrow 0 
\end{equation}
which is compatible with the simplicial structure. Hence we 
have complexes 
\begin{center}
$
\begin{array}{rccccccccc}
({\Pic}^0)^{\bullet}: 
& 
\cdots 
& 
\rightarrow 
&
{\Pic}^{0}(X'_{p-1})
&
{\rightarrow}
&
{\Pic}^{0}(X'_{p})
&
{\rightarrow}
&
{\Pic}^{0}(X'_{p+1}) 
&
\rightarrow
&
\cdots 
\vspace{0.1cm}
\\
{\Pic}^{\bullet}: 
& \cdots 
&
\rightarrow
&
\Pic (X'_{p-1}) 
&
{\rightarrow}
&  
\Pic (X'_{p})
&
{\rightarrow}
&
\Pic (X'_{p+1}) 
&
\rightarrow 
&
\cdots 
\vspace{0.1cm}
\\
{\NS}^{\bullet}: 
& 
\cdots 
&
\rightarrow
& 
{\NS}(X'_{p-1}) 
& 
{\rightarrow}
&
{\NS}(X'_{p})
&
{\rightarrow}
&
{\NS}(X'_{p+1})
& 
\rightarrow 
&\cdots 
\end{array}
$
\end{center}
where the maps are the alternating sums of the maps induced by
the face maps $\delta^i_p: X'_{p+1}\rightarrow X'_{p}, \ 0\leq i 
\leq p+1$. If  $\H^p(-)$ denotes the cohomology of any of the
above complexes, there is a boundary map 
$$\lambda^p: \H^{p-1}(\NS^{\bullet}) \rightarrow 
\H^1({({\Pic^0})}^{\bullet}).$$
The quotient $\H^p({({\Pic^{0}})}^{\bullet})^0$ of the connected 
component of the identity of the group of $p$-cycles of 
${({\Pic^0})}^{\bullet}$, viewed as an algebraic group scheme over 
$k$, by $p$-boundaries is a subgroup of $\H^p({({\Pic^0})}^{\bullet})$. 
The pullback square
$$
\minCDarrowwidth7pt
\begin{CD} 
\L @>>> \H^1({({\Pic^0})}^{\bullet})^0\\
@VVV @VVV\\
\H^{0}(\NS^{\bullet}) @>>> \H^p({({\Pic^0})}^{\bullet})
\end{CD}
$$
defines a finitely generated abelian group $\L$, and an effective 1-motive
\begin{equation}
\Lb(\Xs):=[\L \overset{\lambda}{\rightarrow} 
\H^p({({\Pic^0})}^{\bullet})^0], 
\end{equation}
\cite[\S 1]{BRS}
which we will refer to as the Lefschetz motive. 

\begin{remark} Analogously to $\Lb(\Xs)$ one defines 1-motives 
${\Lb}^p(\Xs):=[\L^{p-1}\rightarrow 
\H^p({(\Pic^0)}^{\bullet})^0]$ 
for $p\geq 1$. The 1-motives $\Lb^p(\Xs)$ have been investigated by the 
first author in \cite{Bash}.  
\end{remark}

\noindent
If ${\cF}_{\d}$ is a simplicial sheaf on $\Xs$ with respect to a topology 
$\tau$, we write ${\E}^{p,q}_r(({\cF}_{\d})_{\tau})$ for the
${\E}^{p,q}_r$-term of the spectral sequence computing the cohomology of 
${\cF}_{\d}$ from the cohomology of its components
\begin{equation}\label{compss}
{\E}^{p,q}_1={\H}^q_{\tau}(X_p, {\cF}_p)
\Rightarrow {\Hb}^{p+q}_{\tau}(\Xs,{\cF}_{\d}).
\end{equation}
Following our convention, if we do not specify $\tau$ it is understood that
we are considering the analytic topology; in some cases when there is no risk 
of confusion, we will also suppress the subscript $\tau$.

\begin{lemma}
\label{motisolemma}
Let $\Xs$ be a smooth proper simplicial scheme over $\C$. 
\begin{itemize}
\item[(i)] ${\E}^{p,q}_2(\Q_{\d})\cong 
\gr^{\W}_q\Hb^{p+q}(\Xs, \Q_{\d}).$ 
\vspace{0.05cm}
\item[(ii)] ${\E}^{p,0}_2(({\Z_{\d})}_{\zar})\cong
\Hb^{p}_{\zar}(\Xs, \Z_{\d}).$ 
Hence $\W_0\hspace{-0.07cm}\Hb^p(\Xs,\Q_{\d})
\cong \Hb^p_{\zar}(\Xs, \Q_{\d})$. 
\end{itemize}
\end{lemma}

\begin{proof} The spectral sequence (\ref{compss}) for the sheaf $\Q_{\d}$
on the analytic site is a spectral sequence of mixed Hodge structures. By
our assumption each component $X_p$ is smooth and proper, thus the groups 
$\H^q(X_p,\Q)$ carry a pure Hodge structure of weight $q$. This implies (i). 
The corresponding spectral sequence for the constant Zariski sheaf $\Z_{\d}$ 
degenerates at $E_2$, which immediately gives (ii). 
\qed \end{proof}

From Lemma \ref{motisolemma} it is easy to see that source and target of the 
1-motives $\Lb(\Xs)$ and $\Lb(\Eb)$ are isomorphic up to torsion. 
We need to show that under these isomorphisms the maps $\e$ and 
$\lambda$ coincide.
\\
\indent 
By \cite[Theorem 3.4]{Beab} the derived category ${\cD}^b(\MHS)$ of bounded
complexes of mixed Hodge structures is equivalent to the triangulated category 
of bounded mixed Hodge complexes \cite[Definition 3.2]{Beab}. 
\\
\indent
Let $\K \rightarrow \K'$ be a morphism of 
bounded mixed Hodge complexes, where $\H^i(\K)$ and $\H^i(\K')$ are pure 
of weight $i$. Let $\Ker:=\ker\{\H^2(\K)\rightarrow \H^2(\K')\}$ and 
$\Cok:=\coker\{\H^1(\K)\rightarrow \H^1(\K')\}$. The distinguished triangle
$$ \rightarrow C(K\rightarrow K') \overset{+1}{\rightarrow} K \rightarrow 
K'\rightarrow $$ 
shows that an element $u\in \Hom_{\MHS}(\Z, \Ker(1))$ gives rise to an 
extension class $\e_u\in \Ext^1_{\MHS}(\Z, \Cok(1)).$
Explicitly $\e_u$ is given by the extension
$$0\rightarrow \Cok(1)
\rightarrow \H^1(\operatorname{C}(\K\rightarrow \K')) \rightarrow \Ker(1) 
\rightarrow 0.$$
Let $u': \Z\rightarrow \K(1)[2]$ be a lift of $u$, and let $u'': \Z\rightarrow 
\K'(1)[2]$ be the composition of $u$ with $\K\rightarrow \K'$. Then 
$u''$ induces an element of $\Ext^1_{\MHS}(\Z, \Cok(1))$ which coincides 
with $\e_u$ \cite[Remark 5.5]{BRS}.  
\\
\indent
Consider the case when $\K=R\Gamma(X_0, \Z)$ and $\K'=R\Gamma(X_1,\Z)$ are the 
Hodge complexes of $X_0$ and $X_1$,~i.e.~the underlying complex of abelian 
groups is the usual chain complex of topological spaces with the analytic 
topology 
\cite[Section 4.3]{Beab}.
Taking the difference of the face maps 
$\delta^0_1, \ \delta^1_1: X_1\rightarrow X_0$ 
induces a map
$$\delta^*_0: R\Gamma(X_0,\Z)
\rightarrow R\Gamma(X_1,\Z).$$
Here $\H^i(R\Gamma(X_0,\Z)$ are the usual singular cohomology groups and, 
by definition, 
$${\H}^i_D(X_p,\Z(k)):=\Hom_{{\cD}^{b}(\MHS)}(\Z, R\Gamma(X_p,\Z)(k)[j])$$
are the Deligne-Beilinson cohomology groups of $X_p$. From Lemma \ref{motisolemma} 
we have
$$\Hb^2(\Xs)^{(1,1)}_{\Q}=\Hom_{\MHS}(\Q, \ker\{\H^2(\K)\rightarrow \H^2(\K')\}(1))
,$$ $$\J^1({\gr}^{\W}_1\Hb^2(\Xs))\subseteq 
\Ext^1_{\MHS}(\Z, \coker\{\H^1(\K)\rightarrow \H^1(\K')\}(1)).$$ 
Now ${\H}^2_D(X_p,\Z(1))\cong \Pic(X_p)$ \cite[Section 5.4]{Beab}, and the map
$\delta^*_0: \Pic(X_0)\rightarrow \Pic(X_1)$ is the map obtained by composition
in the derived category as above. By definition of $\lambda$, if an element $u\in 
\H^0({\NS}^{\bullet})$ lifts to $\Pic(X_0)$, then $\delta^*_0(u)\in \Pic^0(X_1)$
modulo the image of $\Pic^0(X_0)$. Thus the above description of the extension 
class map shows the required compatibility. In summary, we have obtained the 
following.

\begin{proposition}\label{motiso}
Let $\Xs$ be a smooth proper simplicial scheme over $\C$.
Then up to isogeny $\Lb(\Xs)\cong \Lb(\Eb).$ In particular, the 
Hodge realization of the isogeny class of $\Lb(\Xs)$ is $\Eb_\Q$. 
\end{proposition}

We need a different description of the Hodge realization of the 1-motive $\Lb(\Xs)$. 
Let $\Xs$ be any simplicial scheme and let $\omega_{\d}$ be the change of topology 
map from the analytic to the Zariski site of $\Xs$. Consider the simplicial Zariski 
sheaf
$${\cH}^{s}_{\d}(\Z(n)):=\Rb^s(\omega_{\d})_{*}\Z_{\d}(n).$$
The $p$-th component of this sheaf is the Zariski sheaf 
${\cH}^s_{p}(\Z(n))$ associated to the presheaf 
$U\mapsto \H^s(U, \Z(n))$, where $U\subseteq X_p$. For example, 
${\cH}^{0}_{\d}(\Z(n))$ is the constant Zariski sheaf $\Z_{\d}(n)$. The
component spectral sequence (\ref{compss}) is
\begin{equation}\label{compssH1}
{\E}^{p,q}_1={\H}^q_{\zar}(X_p, {\cH}^s_{p}(\Z(n)))\Rightarrow {\Hb}^{p+q}_{\zar}
(\Xs, 
{\cH}^{s}_{\d}(\Z(n))). 
\end{equation}

The Leray spectral sequence with respect to $\omega_{\d}$ has the form
\begin{equation}\label{leray}
{\L}^{p,q}_2=\Hb^p_{\zar}(\Xs, {\cH}^q_{\d}(\Z(n)))\Rightarrow 
\Hb^{p+q}(\Xs,\Z_{\d}(n)).
\end{equation}

There are versions of (\ref{compssH1}) and (\ref{leray}) for rational, 
real and complex coefficients; if there is no risk of confusion, we 
will omit the coefficients from the notation.

\begin{lemma}
\label{H1H1incl}
Let $\Xs$ be a smooth proper simplicial scheme over $\C$. The
Leray spectral sequence (\ref{leray}) induces an exact sequence
$$0\rightarrow \Hb^1_{\zar}(\Xs,{\cH}^1_{\d}(\Q))\rightarrow 
\frac
{
{\Hb}^2(\Xs, \Q_{\d})}
{\W_0 \hspace{-0.06cm}{\Hb}^2(\Xs,\Q_{\d})}
\rightarrow 
\Hb^0_{\zar}(\Xs, {\cH}^2_{\d}(\Q)). $$
\end{lemma}

\begin{proof} Let $\L^{p,q}_{\infty}$ be the filtration induced by the spectral
sequence. We show the differentials $\di_2^{i,1}$ are trivial for $i\geq 0$. 
In particular, $\L^{1,1}_2=\L^{1,1}_{\infty}$, \ $\L_2^{2,0}=\L_{\infty}^{2,0},$
and the exact sequence claimed in the lemma is of the form
$0\rightarrow \L^{1,1}_{\infty} 
\rightarrow {\Hb}^2(\Xs)/\L^{2,0}_{\infty}\rightarrow \L^{0,2}_2$.
Note that the edge map 
$$t^i: \L^{i,0}_2=\Hb^i_{\zar}(\Xs,{\cH}^0_{\d}(\Q))  
\rightarrow \Hb^i(\Xs, \Q_{\d})$$
coincides with the pullback map $\Hb^i_{\zar}(\Xs, \Q_{\d})\rightarrow
\Hb^i(\Xs, \Q_{\d})$. The component spectral sequences of ${({\Q_{\d}})_{\zar}}$ 
and $\Q_{\d}$ induce the commutative diagram
$$
\minCDarrowwidth7pt
\begin{CD}
{\E}^{i,0}_2({(\Q_{\d})}_{\zar}) @>{\cong}>> {\E}^{i,0}(\Q_{\d})\\
@V{\cong}VV @VVV\\
\Hb^i_{\zar}(\Xs,\Q_{\d}) @>{t^i}>> \Hb^i(\Xs,\Q_{\d})
\end{CD}
$$
By Lemma \ref{motisolemma} we can identify $t^i$ with the inclusion
$\W_0\hspace{-0.06cm}\Hb^i(\Xs)\rightarrow \Hb^i(\Xs)$. Thus 
$\L^{i,0}_2=\L^{i,0}_{\infty}$. Our claim follows since $\L^{3,0}_3=
\coker \di_2^{1,1}$ and $\L^{2,0}_{\infty}=\coker \di_2^{0,1}$.
\qed \end{proof}

\begin{lemma}\label{H1H1mix} 
Let $\Xs$ be a smooth proper simplicial scheme over $\C$. There is an
extension of rational mixed Hodge structures
$$
0\rightarrow {\gr}^{\W}_1\Hb^2(\Xs,\Q_{\d})
\rightarrow {\Hb}^1_{\zar}(\Xs, {\cH}^1_{\d}(\Q))
\rightarrow \ker\{ {\NS}(X_0)_{\Q}
\rightarrow {\NS}(X_1)_{\Q}\} \rightarrow 0
$$
which is compatible with the mixed Hodge structure on $\Hb^2(\Xs,\Q_{\d})$. 
In particular, ${\Hb}^1_{\zar}(\Xs, {\cH}^1_{\d}(\Q))$ carries a mixed Hodge 
structure with the properties: 
\begin{itemize}
\item[] $\W_0 \hspace{-0.06cm}{\Hb}^1_{\zar}(\Xs, {\cH}^1_{\d}(\Q))=0,$
\item[] $\W_1 \hspace{-0.06cm}{\Hb}^1_{\zar}(\Xs, {\cH}^1_{\d}(\Q))=
{\gr}^{\W}_1\Hb^2(\Xs,\Q_{\d})$, 
\item[] $\W_2 \hspace{-0.06cm}{\Hb}^1_{\zar}(\Xs, {\cH}^1_{\d}(\Q))=
{\Hb}^1_{\zar}(\Xs, {\cH}^1_{\d}(\Q))$, 
\item[] ${\gr}^{\W}_2{\Hb}^1_{\zar}(\Xs, {\cH}^1_{\d}(\Q))
=\ker\{ {\NS}(X_0)_{\Q}\rightarrow {\NS}(X_1)_{\Q}\}=
{\gr}^{\W}_2 \Hb^2(\Xs,\Q_{\d})^{(1,1)}.$
\end{itemize}
\end{lemma}

\begin{proof} We claim we have a commutative diagram with exact rows
$$
\minCDarrowwidth7pt
\begin{CD}
0 @>>> {\Hb}^1_{\zar}(\Xs, {\cH}^1_{\d}(\Q)) @>>> 
{\Hb}^2(\Xs,\Q_{\d})/\W_0{\Hb}^2(\Xs,\Q_{\d})
@>>> {\Hb}^{0}_{\zar}(\Xs,{\cH}^2_{\d}(\Q))\\
@. @V{\alpha}VV @V{\beta}VV @VV=V\\
0 @>>>
{\E}^{0,1}_2({\cH}^1_{\d}(\Q))
@>>>
{\E}^{0,2}_2(\Q_{\d})
@>>>
{\E}^{0,2}_2({\cH}^2_{\d}(\Q))
\end{CD}
$$
whose vertical maps are induced by the componentwise spectral sequences, 
and whose horizontal maps are coming from the Leray spectral sequences
along the change of topology maps for ${\Xs}$ and its components $X_p$. 
To see that this diagram commutes, consider first the constant simplicial
scheme $(X_0)_{\d}$. The map $\Xs\rightarrow (X_0)_{\d}$ clearly induces
such a commutative diagram with the bottom row replaced by the
corresponding 
${\E}^{0,i}_1$-terms. Since for any simplicial sheaf $\cF_{\d}$, 
$${\E}^{0,i}_2({\cF}_{\d})=\ker\{\H^i(X_0,{\cF}_0)\rightarrow \H^i(X_1,{\cF}_1)\}
\subseteq {\E}^{0,i}_1({\cF}_{\d})=\H^i(X_0,{\cF}_0),$$
this implies our claim. Now we use that by Lemma \ref{motisolemma} we have the 
formula
$${\E}^{0,2}(\Q_{\d})=\ker\{\H^2(X_0,\Q)\rightarrow 
\H^2(X_1,\Q)\}=
{\gr}^{\W}_2\Hb^2(\Xs,\Q_{\d}).$$
Hence
$$\ker \alpha =\ker \beta ={\gr}^{\W}_1\Hb^2(\Xs,\Q_{\d})$$
and a diagram chase shows ${\Hb}^1_{\zar}(\Xs, {\cH}^1_{\d}(\Q))$ 
surjects onto the group
\begin{equation}\label{gr11piece}
{\E}^{0,1}({\cH}^1_{\d}(\Q))=
\ker\{\NS(X_0)_{\Q}\rightarrow
\NS(X_1)_{\Q}\}={\gr}^{\W}_2\Hb^2(\Xs,\Q_{\d})^{(1,1)}.
\end{equation}
\indent
In particular we have a commutative diagram of mixed Hodge structures
\begin{equation}
\label{MHSdia}
\minCDarrowwidth6pt
\begin{CD}
0 
@>>> 
{\gr}^{\W}_1\Hb^2(\Xs,\Q_{\d})
@>>>
{\Hb}^1_{\zar}(\Xs, {\cH}^1_{\d}(\Q))
@>>>
{\gr}^{\W}_2\Hb^2(\Xs, \Q_{\d})^{(1,1)}
@>>> 
0\\
@. @V{=}VV @VVV @VVV\\
0 @>>> {\gr}^{\W}_1\Hb^2(\Xs,\Q_{\d}) @>>> \Hb^2(\Xs, \Q_{\d})/\W_0
@>>> {\gr}^{\W}_2\Hb^2(\Xs,\Q_{\d}) @>>> 0. 
\end{CD}
\end{equation}

\qed \end{proof}

\begin{proposition}
\label{H1H1ishodge}
Let $\Xs$ be a smooth proper simplicial scheme over $\C$. Then
${\Hb}^1_{\zar}(\Xs, {\cH}^1_{\d}(\Q))$ is isomorphic to the 
Hodge realization of the isogeny class of the Lefschetz 1-motive
$\Lb(\Xs)$. 
\end{proposition}

\begin{proof} By (\ref{Ebdef}), (\ref{gr11piece}) and (\ref{MHSdia}), the 
mixed Hodge structure $\Eb_{\Q}$ is canonically isomorphic to 
${\Hb}^1_{\zar}(\Xs, {\cH}^1_{\d}(\Q))$. Now use Proposition \ref{motiso}. 
\qed \end{proof}


\section{Cohomology classes of simplicial line bundles}
\noindent
Let $\Xs$ be a smooth proper simplicial scheme over an algebraically
closed field $k$. By \cite[4.1]{BaSrone} the $fppf$-sheafification of 
the simplicial Picard functor $T\mapsto \Pic(\Xs\times_k T)$ is representable, 
and its connected component of the identity 
$\Pic^0(\Xs)$ is a semi-abelian variety over $k$. Explicitly, the semi-abelian
variety $\Pic^0(\Xs)$ is given by the extension 
\begin{equation}\label{Picoext}
0\rightarrow \T(\Xs)\rightarrow \Pic^0(\Xs) \rightarrow \A(\Xs)\rightarrow
0,
\end{equation}
where $\A(\Xs)$ is the connected component of $\ker\{\Pic^0(X_0)\rightarrow 
\Pic^0(X_1)\}$, and 
\begin{equation}
\label{T=E2}
\T(\Xs)=\E^{1,0}_2({\OO}^*_{\Xs})=\Hb^1_{\zar}(\Xs,k^*_{\d})
\end{equation}
is the connected $k$-torus which is the $\E^{1,0}_2$-term of the 
spectral sequence computing $\Pic(\Xs)$. 
We define the simplicial N\'eron-Severi group as the quotient 
$$\NS(\Xs):=\coker\{\Pic^0(\Xs)\rightarrow \Pic(\Xs)\}.$$
From the above description of $\Pic^0(\Xs)$ it is easy to see
that $\NS(\Xs)$ injects into $\NS(X_0)$ up to a finite group;
in particular, the simplicial N\'eron-Severi group is a finitely
generated abelian group. 

If $\Xs$ is a smooth proper simplicial scheme over $\C$ we have
two cycle class maps: The simplicial exponential 
sequence induces the analytic cycle map
\begin{equation}
\c_1:\Pic(\Xs)\rightarrow \Hb^2(\Xs, \Z_{\d}(1)),
\end{equation}
and if $\omega_{\d}$ is the change of topology map from the analytic to the
Zariski site, the restriction of 
$(\omega_{\d})_*(\OO_{\Xs}^*)_{\an}
\rightarrow 
{\R}^1(\omega_{\d})_*\Z_{\d}(1)$ to the 
subsheaf $(\OO^*_{\Xs})_{\zar}$ induces
\begin{equation}
\tilde{\c_1}: \Pic(\Xs)\rightarrow \Hb^1_{\zar}(\Xs, {\cH}^1_{\d}(\Z_{\d}(1))).
\end{equation}
Note that the map $\OO^*_{\Xs}\rightarrow {\cH}^1_{\d}(\Z(1))$ defining 
$\tilde{\c}_1$ coincides with the composition
$${\OO}_{\Xs}^*\rightarrow {\cF}^{1,1}_{\d} \rightarrow {\cH}^1_{\d}(\Z(1)),$$
where ${\cF}^{1,1}_{\d}$ is the simplicial sheaf whose $p$-the component is the
Zariski sheaf associated to the presheaf defined by sending an open subset 
$U\subseteq X_p$ to 
$$U \mapsto {\cF}^{1,1}_p(U):=\Hom_{\MHS}(\Z(0), \H^1(U,\Z(1))\subseteq 
\H^1(U,\Z(1)).$$ 
In particular, the map $dlog: \OO^*_{X_p}(U)\rightarrow {\cF}^{1,1}_p(U), \ 
f\mapsto df/f$ induces the following exact sequence
\begin{equation}
\label{F11exact}
0\rightarrow \C^*_{\d} \rightarrow {\OO}_{\Xs}^* \rightarrow {\cF}^{1,1}_{\Xs} 
\rightarrow 0.
\end{equation}
\begin{lemma} 
\label{F11prop}
Let $\Xs$ be a smooth proper simplicial scheme over $\C$. 
\begin{itemize}
\item[(i)] $\Hb^0_{\zar}(\Xs, {\cF}^{1,1}_{\d})=0$
\vspace{0.05cm}
\item[(ii)] $\Hb^{1}_{\zar}(\Xs, {\cF}^{1,1}_{\d})=\ker\{\Pic(X_0)
\rightarrow \Pic(X_{1})\}.$
\vspace{0.05cm}
\item[(iii)] $\Hb^0_{\zar}(\Xs,{\cH}^1_{\d}/{\cF}^{1,1}_{\d})=
\ker\{\H^1(X_0,{\OO}_{X_0})
\rightarrow \H^1(X_1,{\OO}_{X_1})\}$. 
\end{itemize}
\end{lemma}

\begin{proof} For any smooth proper variety  $X_p$ the sheaf ${\cF}^{1,1}_{p}$ has no 
global sections, thus (i). Clearly, $\H^r(X_p,{\cF}^{1,1}_{p})=0$ for $r\geq 2$. Since
$\H^1(X_p,{\cF}^{1,1}_{p})=\Pic(X_p)$ \cite[Theorem 1.3]{Esnote} (ii) is immediate 
from the componentwise spectral sequence. For (iii) use that the quotient sheaf
${\cH}^1_{p}/{\cF}^{1,1}_{p}$
is a constant sheaf with value $\H^1(X_p,{\OO}_{X_p})$, see \cite{BaSrBloch}. 
\qed \end{proof}

\begin{lemma}\label{cycleclasskernel}
Let $\Xs$ be a smooth proper simplicial scheme over $\C$. Then
$$\Pic^0(\Xs)=\ker c_1=\ker \tilde{\c}_1.$$
\end{lemma}

\begin{proof} For $\c_1$ this follows from the simplicial exponential 
sequence and the fact that the Hodge filtration on $\Hb^*(\Xs,\C_{\d})$
is defined by truncations of the simplicial De Rham complex 
\cite[8.1.12]{DIII}; in particular, $\Hb^2(\Xs, \OO_{\Xs})=\Hb^2(\Xs)/\F^1$. 

For $\tilde{\c}_1$ note that
the restriction of $\tilde{\c}_1$ to $\Pic^0(\Xs)$ and the composition 
$\ker \tilde{\c}_1\rightarrow \Pic(\Xs)\rightarrow \NS(\Xs)$ are trivial, 
since $\NS(\Xs)$ and $\im \tilde{\c}_1$ are finitely generated abelian groups,
and $\ker \tilde{\c}_1$ is divisible by Lemma \ref{F11prop}.  
\qed \end{proof}

\begin{lemma}\label{NSmix}
Let $\Xs$ be a smooth proper simplicial scheme over $\C$. 
Then
$$\F^1\Hb^2(\Xs)\cap \Hb^2(\Xs, \Z_{\d}(1))
\cong \NS(\Xs)\subseteq \Hb^1_{\zar}(\Xs, {\cH}^1_{\d}(\Z(1))),$$
where `intersection' means the inverse image under $\Hb^2(\Xs,\Z_{\d})
\rightarrow \Hb^2(\Xs, \C_{\d})$. 
\end{lemma}

\begin{proof} Since
$\Hb^2(\Xs, \OO_{\Xs})\cong \Hb^2(\Xs,\C_{\d})/\F^1$ the exponential 
sequence shows
$$\im\{\Pic(\Xs)\rightarrow \Hb^2(\Xs, \Z_{\d}(1))\}
=
\ker\{\Hb^2(\Xs, \Z_{\d}(1))\rightarrow 
\Hb^2(\Xs, \C_{\d})/\F^1\}.$$
Hence $\NS(\Xs)=\F^1\cap \ \Hb^1(\Xs,\Z_{\d})$. The remaining claim is a 
consequence of the fact that the kernel of both maps $\c_1$ and $\tilde{\c}_1$ 
is isomorphic to $\Pic^0(\Xs)$. 
\qed \end{proof}

Define 
$$\Ib(\Xs):=\coker\{\NS(\Xs)\rightarrow 
\Hb^1_{\zar}(\Xs, {\cH}^1_{\d}(\Z))\}.$$
Note that $\Ib(\Xs)$ is the quotient of a finitely generated abelian group
and $\Ib(\Xs)_{\Q}$ inherits a mixed Hodge structure from 
$\Hb^1_{\zar}(\Xs, {\cH}^1_{\d}(\Q))$. 

\begin{lemma}\label{Imix}
Let $\Xs$ be a smooth proper simplicial scheme over $\C$. 
\begin{itemize}
\item[] $\W_0\Ib(\Xs)_{\Q}=0,$
\vspace{0.05cm}
\item[] $\W_1\Ib(\Xs)_{\Q}={\gr}^{\W}_1\Hb^2(\Xs,\Q_{\d})$, 
\vspace{0.05cm}
\item[] $\W_2\Ib(\Xs)_{\Q}=\Ib(\Xs)_{\Q}$, 
\vspace{0.05cm}
\item[] ${\gr}^{\W}_2\Ib(\Xs)_{\Q}=
{\gr}^{\W}_2\Hb^2(\Xs,\Q_{\d})^{(1,1)}/(\F^1\cap
(\Hb^2(\Xs,\Q_{\d})/\W_0))=\im \lambda_{\Q}$. 
\end{itemize}
\end{lemma}

\begin{proof} The functor $\W_p(-)$ is exact, thus all claims about
the weight subspaces follow from Lemma \ref{H1H1mix}. For the last
claim we use the isomorphisms
$${\gr}^{\W}_2\Ib(\Xs)_{\Q}\cong
\frac
{{\gr}^{\W}_2\Hb^1_{\zar}(\Xs {\cH}^1_{\d}(\Q))}
{
F^1\cap \  \Hb^1_{\zar}(\Xs, {\cH}_{\d}^1(\Q))}
\cong
\frac
{{\gr}^{\W}_2\Hb^2(\Xs, \Q_{\d})^{(1,1)}}
{\F^1\cap \ (\Hb^2(\Xs, {\Q_{\d}})/\W_0)}
\cong
\im \e_{\Q}.$$
Here the first and the last identification are clear, and the second one
follows from Lemma \ref{motisolemma}, Lemma \ref{H1H1mix} and Proposition 
\ref{motiso}. In particular, we note
\begin{equation}\label{kerneliso}
\ker \lambda_{\Q} \cong \F^1\cap \  \Hb^1_{\zar}(\Xs, {\cH}_{\d}^1(\Q)) 
\cong 
{\F^1\cap \ (\Hb^2(\Xs, \Q_{\d})/\W_0)}.
\end{equation}
\qed \end{proof}


\section{$(1,1)$-classes on proper seminormal varieties}
\label{11classes}
\noindent
We first give a characterization of seminormality in terms of hypercoverings;
we believe this is known to experts although we could not find a reference in
the literature. 

\begin{lemma}\label{seminormal-simplicial} 
Let $X$ be a proper reduced scheme over $\C$, and let
$\pi_{\d}: \Xs \rightarrow X$ be a smooth proper hypercovering. The following
statements are equivalent.
\begin{itemize}
\item[(i)] $X$ is seminormal.
\vspace{0.05cm}
\item[(ii)] 
$\pi_{\d*}{\OO}_{\Xs}\cong {\OO}_{X}.$
\end{itemize}
\end{lemma}

\begin{proof} Let $X_{\sn}$ be the seminormalization of $X$. Then $X_{\sn}$ satisfies
the universal property that any morphism $Y\rightarrow X$ with $Y$ seminormal factors
uniquely through $X_{\sn}$ \cite[Chapter I, Proposition 7.2.3.3]{Kollar}. 
Since $X_{\sn}\rightarrow X$ is bijective, the functor from
$X_{\sn}$-schemes to $X$-schemes is faithful. The structure
maps $X_i\rightarrow X$ factor through $X_0$, thus by the universal property
of seminormalization through $X_{\sn}$. Furthermore, the maps $X_i\rightarrow
X_{\sn}$ are compatible with the simplicial structure,~i.e.~$\Xs$ is a 
simplicial $X_{\sn}$-scheme in a canonical way, and $\pi_{\d}$ factors as 
$$
\Xs\rightarrow X_{\sn} \rightarrow X.
$$
Let $\widetilde{X}=\Spec(\pi_{\d *}\OO_{\Xs})$, and consider the factorization of
$\pi_{\d}:\Xs \rightarrow X$ as
$$
\Xs \rightarrow \widetilde{X}\rightarrow X.
$$
We claim that $X_{\sn}=\widetilde{X}$. Note that $\widetilde{X}\rightarrow X$ clearly
factors through $X_{\sn}$ (with the obvious notation, $\widetilde{X_{\ns}}=\widetilde{X}$). 
On the other hand, if $(X_0)_{\an}/(X_1)_{\an}$ denotes the topological
quotient given by the equivalence relation determined by $(X_1)_{\an}$, we have
continuous proper maps
$$X_{\an}\leftarrow (X_0)_{\an}/(X_1)_{\an} \rightarrow (\tilde{X})_{\an}.$$
Here the first map is a bijection (and thus a homeomorphism) by definition of 
a hypercovering, and 
the second map results from the fact that $X_0$ is an $\tilde{X}$-scheme.
The continuous section $X_{\an}\rightarrow (\tilde{X})_{\an}$ implies that 
the finite birational morphism $\tilde{X}\rightarrow X$ must be a bijection
on closed points, and therefore $X_{\sn}\rightarrow X$ factors through $\widetilde{X}
\rightarrow X$. 
\qed \end{proof}

Assume $X$ is a proper reduced scheme over $\C$, and  
$\pi_{\d}: \Xs\rightarrow X$ is a smooth proper hypercovering.
We consider divisor classes on $X$ and $\Xs$. By \cite[Remark 4.11]{BaSrone}
$\Pic^0(\Xs)$ is independent of the choice of $\pi_{\d}$.
Lemma \ref{NSmix} shows that the 
same holds for $\NS(\Xs)$. Hence $\Pic(\Xs)$ is independent of the choice of 
a hypercovering. Lemma \ref{Imix} and its proof show that the 
1-motive $\Lb(\Xs)$ as well is independent of the choice of $\pi_{\d}$.

\begin{lemma} 
\label{exponentialseq}
Let $X$ and $\pi_{\d}: \Xs\rightarrow X$ be as above. Then
\begin{itemize}
\item[(i)] $\NS(X)\subseteq \NS(\Xs)$ and 
$\NS(\Xs)/\NS(X)\cong \Pic(\Xs)/\Pic(X)$ 
which is a finite lattice (~i.e.~isomorphic to $\Z^k$ for some integer $k$).

\vspace{0.05cm}
\item[(ii)] 
If $\Pic(\Xs)$ is generated by line bundles on 
$X$, $\NS(X)=\F^1\cap \H^2(X,\Z)$.

\vspace{0.05cm}
\item[(iii)] $\pi_{\d}^*:\Pic^0(X)\rightarrow \Pic^0(\Xs)$ is surjective, with kernel
a $\C$-vector space. 
\vspace{0.05cm}
\item[(iv)] 
If $\H^p(X,\OO_X)=\H^p(X,\C)/\F^1$ for $p=1,2$, then $\Pic(X)\cong \Pic(\Xs)$. 
\vspace{0.05cm}
\item[(v)] If $X$ is seminormal, 
$\Pic^0(X)\cong \Pic^0(\Xs)\cong \H^1(X,\C)/({\F}^1+\H^1(X,\Z)).$ 
\end{itemize}
\end{lemma}

\begin{proof} Consider the map $\H^p(X,\C)\rightarrow \H^p(X,\OO_X)$ induced 
by the inclusion $\C\rightarrow \OO_X$. The isomorphism 
$\Hb^p(\Xs,{\OO}_{\Xs})\cong 
{\gr}^0_{\F}\  \Hb^p(\Xs, \C_{\d})$ and descent imply the map 
$\pi^*_{\d}:\H^p(X, {\OO}_X)\rightarrow 
{\Hb}^p(\Xs,{\OO}_{\Xs})$ is surjective 
for all $p\geq 0$. Comparing the exponential sequences for $X$ and $\Xs$ using 
GAGA we a obtain a commutative diagram with exact rows
$$
\minCDarrowwidth6pt
\begin{CD}
@>>> \H^1(X,\Z) @>>> \H^1(X,{\OO}_X) @>>> \Pic(X) @>>> \H^2(X,\Z) 
@>>> \H^2(X,{\OO}_X) @>>> \\ 
@. @V{\cong}VV @V{\text{onto}}VV @VVV @VV{\cong}V @VV{\text{onto}}V\\
@>>> \Hb^1(\Xs,\Z_{\d}) @>>> \Hb^1(\Xs,{\OO}_{\Xs}) @>>> \Pic(\Xs) 
@>>> \Hb^2(\Xs,\Z_{\d}) @>>> \Hb^2(\Xs,{\OO}_{\Xs}) @>>>
\end{CD}
$$
Now (i)-(iv) follow from a diagram chase using descent and Lemma \ref{NSmix}. 
For (v) note that if  $X$ is seminormal we have an isomorphism 
$\pi_{\d*}{\OO}^*_{\Xs}\cong {\OO}^*_{X}$ by Lemma
\ref{seminormal-simplicial}. 
The Leray spectral sequence 
of $\pi_{\d}$ implies $\pi^*_{\d}: \Pic(X)\rightarrow \Pic(\Xs)$ is
injective; combined with (iii) this shows $\Pic^0(X)\cong \Pic^0(\Xs)$. 
\qed \end{proof}

Let $X$ be a proper reduced scheme over $\C$ with 
smooth proper hypercovering $\pi_{\d}: \Xs \rightarrow X$. 
Consider the commutative diagram with exact rows
\begin{equation}
\label{H2locdia}
\minCDarrowwidth7pt
\begin{CD}
   0@>>> \NS(X) @>>> \H^2_{\zarloc}(X, \Z) @>>> \I(X)_0 @>>> 0\\
   @. @V{\pi^*_{\d}}VV @V{\pi^*_{\d}}VV @VV{\pi^*_{\d}}V\\
   0@>>> \NS(\Xs) @>>> \Hb^1_{\zar}(\Xs, {\cH}^1_{\d}) @>>>
   \Ib(\Xs) @>>> 0\
\end{CD}
\end{equation}

Here the top row is the exact sequence constructed in \cite{BaSr1}
(modified to allow for $X$ being reducible). The bottom row is obtained
from Lemma~\ref{NSmix}. The quotient $\I(X)_0$ is known to be a lattice 
of finite rank. The above diagram is a diagram of mixed Hodge structures. 
From the snake lemma we have an induced map
\begin{equation}
\label{deltamap}
\delta: \frac{\NS(\Xs)}{\NS(X)} 
\rightarrow 
\frac{
\Hb^1_{\zar}(\Xs, {\cH}^1_{\d})}
{\im \H^2_{\zarloc}(X, \Z)}.
\end{equation}

\begin{lemma}\label{deltalemma}
Let $X$ be a proper reduced scheme over $\C$. Then 
$$\NS(X)_{\Q}={\F}^1 \cap \H^2_{\zarloc}(X,\Q) 
\Leftrightarrow \ker {\delta}_{\Q}=0.$$
\end{lemma}

\begin{remark}\label{torsion}
Since $X$ is proper, it follows from the exponential sequence and GAGA that 
the torsion in $\H^2(X,\Z)$ is always algebraic. The above lemma implies that 
$\NS(X)={\F}^1 \cap \H^2_{\zarloc}(X,\Z)$, provided $\ker {\delta}_{\Q}=0$. 
\end{remark}

\begin{proof} 
By descent and Lemma \ref{H1H1incl}, the kernel of
$\H^2_{\zarloc}(X,\Q)\rightarrow \Hb^1_{\zar}(\Xs, {\cH}^1_{\d}(\Q))$ 
in (\ref{H2locdia}) is precisely 
${\W}_0\hspace{-0.1cm}\H^2_{\zarloc}(X)$. Lemma \ref{NSmix} and 
Lemma \ref{Imix} show that
$$\ker {\delta}_{\Q}=
(\F^1\cap (\H^2_{\zarloc}(X,\Q)/{\W_0}))/\NS(X)_{\Q}.$$ 
The claim follows since ${\F}^1\cap {\W}_0=0$. 
\qed \end{proof}

To prepare for the proof of our main theorem we need three lemmas.

\begin{lemma}\label{restriction} Let $\pi_{\d}: \Xs \rightarrow X$ be
a smooth proper hypercovering. Let $A$ denote either $\Z, \Q, \C$ or 
$\C^*$. Then the restriction map 
$$\Hb^p_{\zar}(\Xs,A_{\d})_X\rightarrow {\Rb}^p\pi_{\d *}A_{\d}$$ 
is surjective for $p\geq 0$. In particular, the sheaves 
${\Rb}^p\pi_{\d *}\Z_{\d}$ are $\Z$-constructible in the Zariski topology 
(see \cite{BaSr1} for the definition).   
\end{lemma}

\begin{proof} We will only consider $A=\C^{*}$; the other cases are similar. 
If $U\subseteq X$ is Zariski open, let $U_{\d}=\pi^{-1}_{\d}(U)\subseteq \Xs$ 
be the inverse image,  and let $Y_{\d}=\Xs \setminus U_{\d}$ be the complement.
We need to show the following map (of Zariski cohomology groups) is injective
\begin{equation}\label{injsupport}
\Hb^{p+1}_{Y_{\d}}(\Xs,\C^*_{\d})\rightarrow \Hb^{p+1}_{\zar}(\Xs, \C^*_{\d}). 
\end{equation}
In the spectral sequences computing these groups $\E_1^{pq}=0$ for $q\geq 1$.
Hence 
\begin{itemize}
\item[] \hspace{3cm} $\Hb^{p+1}_{Y_{\d}}(\Xs,\C^*)=\H^{p+1}(\H^0_{Y_{\d}}
(\Xs,\C^*_{\d})^{\bullet}),$
\vspace{0.05cm}
\item[] \hspace{3cm} $\Hb^{p+1}_{\zar}(\Xs,\C^*)=\H^{p+1}(\H^0
(\Xs,\C^*_{\d})^{\bullet}),$
\end{itemize}
where $\H^0_{Y_{\d}}(\Xs,\C^*_{\d})^{\bullet}$ and $\H^0_{\zar}
(\Xs,\C^*_{\d})^{\bullet}$ are the complexes of global sections induced by 
the simplicial structure. Consider the map of complexes 
\begin{equation}\label{simplcomplex}
\H^0_{Y_{\d}}(\Xs,\C^*_{\d})^{\bullet}\rightarrow 
\H^0(\Xs,\C^*_{\d})^{\bullet}
\end{equation}
\noindent
We will show (\ref{simplcomplex}) is split injective, this implies (\ref{injsupport}) 
is injective as claimed. 
\\
\indent
Note that for each $n$ the group $\H^0(X_n,\C^*)$ is the free abelian group 
on connected components of $X_n$, tensored with $\C^*$. Similarly, 
$\H^0_{Y_n}(X_n,\C^*)$ is the free abelian group on connected components of 
$X_n$ which are contained in $Y_n$, tensored with $\C^*$. There is an obvious 
inclusion between these groups, and a natural splitting. We need to show this 
splitting is compatible with the differentials. That is, given a $\C^*$-term 
in the kernel of the splitting at level $n$ (corresponding to a component $Z$
of $X_n$ which is not supported in $Y_n$), the face maps send this term into 
a $\C^*$-term in the kernel of the splitting at level $n-1$ (corresponding to 
a component $W$ of $X_{n-1}$ which is not supported in $Y_{n-1}$). 
\\
\indent
Recall that $Y_n=X_n\setminus U_n$, where $U_n$ is the inverse image of $U$ under
the structure morphism $\pi_n: X_n\rightarrow X$. Consider a $\C^*$-term in 
$\H^0(X_n,\C^*)$ that corresponds to a component $Z$ of $X_n$ which is not 
supported in $Y_n$. Thus $Z\cap U_n\neq \emptyset$ and $\pi_n(Z)\cap U\neq 
\emptyset$. Assume $W$ is a component of $X_{n-1}$ such that there is some face map 
$Z\rightarrow W$.  Then this face map is in fact a morphism of $X$-schemes, and the 
image $\pi_{n-1}(W)\subset X$ contains the image $\pi_n(Z)\subset X$. In particular, 
$W$ is {\it not} supported in $Y_{n-1}$, and the $\C^*$-term in $\H^0(X_{n-1},\C^*)$ 
corresponding to $W$ lies in the kernel of the splitting at level $n-1$. 
\qed \end{proof}

\begin{lemma} 
\label{R1inj}
Let $\pi_{\d}:\Xs\rightarrow X$ be a smooth proper hypercovering. The
map of simplicial sheaves ${\OO}^*_{\Xs}\rightarrow {\cF}^{1,1}_{\d}$
induces an injective map $\Rb^1\pi_{\d *}{\OO}^*_{\Xs}\rightarrow 
\Rb^1\pi_{\d *}{\cF}^{1,1}_{\d}$.  
\end{lemma}

\begin{proof} If $A$ is an abelian group, we write $A_X$ for the constant 
sheaf $A$ on $X$. From the description of $\Pic^0(\Xs)$ in (\ref{Picoext}) 
and (\ref{T=E2}) we have an injective map $\Hb^1_{\zar}(\Xs,\C^*_{\d})
\rightarrow \Pic^0(\Xs)$ which fits into the commutative diagram
$$\minCDarrowwidth7pt
\begin{CD}\Hb^1_{\zar}(\Xs, \C^*_{\d})_X @>>> \Pic^0(\Xs)_{X}\\
@VVV @VVV \\
\Rb^1\pi_{\d *}\C^*_{\d} @>>> \Rb^1\pi_{\d *}{\OO}^*_{\Xs} @>>> 
\Rb^1\pi_{\d *}{\cF}^{1,1}_{\d}\\
\end{CD}
$$
whose vertical maps are the sheafification maps. The bottom row is induced by 
(\ref{F11exact}) and is exact. We claim $\Pic^0(\Xs)\rightarrow 
{\Rb}^1{\pi}_{\d *}{\OO}^*_{\Xs}$ is the zero map. Let ${\OO}_{X,x}$ be the local 
ring of a point $x\in X$. Consider the commutative diagram
$$
\minCDarrowwidth6pt
\begin{CD}
\Pic^0(\Xs)@>>> \Pic(\Xs) @>>> \H^0(X,\Rb^1\pi_{\d *}{\OO}_{\Xs})\\
@. @VVV @VVV\\
@. \Pic(\pi^{-1}_{\d}(\Spec({\OO}_{X,x})) @>>> 
(\Rb^1\pi_{\d *}{\OO}_{\Xs})_{x}
\end{CD}
$$
By Lemma \ref{exponentialseq}(iii) $\pi^*_{\d}: \Pic^0(X)\rightarrow 
\Pic^0(\Xs)$ is surjective. Thus any element of $\Pic^0(\Xs)$ lifts to 
an element of $\Pic(X)$ and has trivial image in 
$\H^0(X,{\Rb}^1\pi_{\d *}{\OO}_{\Xs}^*)$; this shows the map in question is 
trivial on stalks. 
\\
\indent
To finish the proof, note that by the previous lemma 
$\Hb^1_{\zar}(\Xs,\C^*_{\d})_X \rightarrow {\Rb}^1\pi_{\d *}\C^*_{\d}$ is 
surjective, thus  
$\Rb^1\pi_{\d *}{\OO}^*_{\Xs}\rightarrow 
\Rb^1\pi_{\d *}{\cF}^{1,1}_{\d}$ is injective. 
\qed \end{proof}

\begin{lemma}\label{piclemma}
Let $f:Y\to Z$ be a morphism between $k$-varieties, where $k$ is an 
algebraically closed field, and $Y$ is non-singular and proper over 
$k$. Let $z\in Z$ be a (closed) point, $F=f^{-1}(z)$ the reduced fiber,
$\tilde{F}$ the seminormalization of $F$. Let 
$Y_z= {\rm Spec}\,({\mathcal O}_{Z,z})\times_ZY$. 
Consider the commutative triangle of groups
$$
\xymatrix{
{\Pic}^0(Y) 
{\ar[r]^{\alpha}}  
{\ar[dr]_{\beta}}
& {\ar[d]} {\Pic(Y_z)}\\
& 
{\Pic(\widetilde{\F})} 
}
$$
Then $\alpha$ and $\beta$ have the same kernels, and isomorphic images. 
Further, the image of $\beta$ is contained in ${\Pic}^0(\widetilde{F})$. 
\end{lemma}
\begin{proof} 
It clearly suffices to prove that $\alpha$ and $\beta$ have the same
kernels. Let $\widehat{{\mathcal O}_{z,Z}}$ denote the completion of 
${\mathcal O}_{z,Z}$ and set 
$\widehat{Y_z}={\rm Spec}\,\widehat{{\mathcal O}_{z,Z}}\times_ZY.$
Since cohomology commutes with flat base change, we see at once that
\[{\rm Pic}(Y_z)\to {\rm Pic}(\widehat{Y_z})\]  
is injective. If ${\mathfrak M}_z\subseteq {\mathcal O}_{z,Z}$ is the 
maximal ideal, let 
$$F_n={\rm Spec}\frac{{\mathcal O}_{z,Z}}{{\mathfrak M}_z^n}\times_ZY.$$
The Formal Function Theorem \cite[III, Theorem 11.1]{Ha} implies the 
natural map
$${\rm Pic}(\widehat{Y_z})\to \liminv{n}{\rm Pic}(F_n)$$
is injective (in fact this map is an isomorphism by \cite[II, Ex.9.6]{Ha}). 
\\
\indent
It suffices to show that, for any $n\geq 1$, the maps of abelian groups
\[{\rm Pic}^0(Y)\to {\rm Pic}(F_n)\]
and 
\[{\rm Pic}^0(Y)\to {\rm Pic}(\widetilde{F})\]
have the same kernel. The kernel of the second map is the group of
$k$-points of a reduced group subscheme of ${\rm Pic}^0(Y)$, which is
hence an extension of a finite group by an abelian variety. 
\\
\indent
In fact, since $Y$, $F_n$ and $\widetilde{F}$ are all proper over $k$, the
respective Picard groups are $k$-points of group schemes over $k$, and
the homomorphisms between Picard groups correspond to homomorphisms of
group schemes. Here ${\rm Pic}^0(Y)$ is in fact an abelian
variety. 
\\
\indent
The morphism $\widetilde{F}\to F_n$ factorizes as
\[\widetilde{F}\to (F_1)_{\rm red}\to F_1\to F_n,\]
where the first arrow is the seminormalization map \cite{Sw}, and the 
other two are infinitesimal extensions (bijective closed immersions defined 
by nilpotent ideal sheaves). The kernels of the induced homomorphisms on 
Picard schemes
\[{\rm Pic}(F_n)\to {\rm Pic}(F_1)\to {\rm Pic}((F_1)_{\rm red})\to {\rm
Pic}(\widetilde{F})\]
are vector group schemes (with affine spaces over $k$ as underlying scheme).
This holds for $(F_1)_{\rm red}\to F_n$ by \cite[VI, Proposition 4.15]{Bo}.
For $\widetilde{F}\to (F_1)_{\rm red}$ we use that by Lemma 
\ref{exponentialseq}(iii) (note that the map $\Xs \to X$ factors through 
the seminormalization),
the map ${\rm Pic}((F_1)_{\rm red})\to {\rm Pic}(\widetilde{F})$ is an 
isomorphism on torsion, thus the kernel is an affine group scheme which
is torsion free,~i.e.~a vector group scheme. 
\\
\indent
The lemma now follows from the fact that if $G$ is an extension of a
finite group by an abelian variety over $k$, and $H$ is a vector group
scheme over $k$, then any homomorphism of groups schemes $G\to H$ is
trivial. 
\qed \end{proof}

Let $\pi_{\d}:\Xs\rightarrow X$ be a smooth proper hypercovering. 
If $\cF_{\d}$ is a sheaf on $\Xs$, we write  
${\E}^{p,q}_2(\Rb^1\pi_{\d *}{\cF}_{\d})$ for the sheaf corresponding 
to the $\E_2^{p,q}$-term of the componentwise spectral sequence 
computing the sheaf $\Rb^1\pi_{\d *}{\cF}_{\d}$. For example, 
$${\E}^{0,1}_2(\Rb^1\pi_{\d *}{\OO}_{\Xs}^*)=
\ker\{\R^1\pi_{0*}{\OO}^*_{X_0} \rightarrow 
\R^1\pi_{1*}{\OO}^*_{X_1}\}.$$

We are ready to prove our main theorem.
\begin{proof}(of Theorem \ref{seminormallefschetz})
\vspace{0.1cm}
\\
By Lemma \ref{deltalemma} and Remark \ref{torsion} it suffices to show 
$\ker \delta=0$. Our first step is to use a Leray spectral sequence argument 
to produce an {\it injective} map
\begin{equation}\label{NSquotientinj}
\NS(\Xs)/\NS(X) \rightarrow 
\H^0(X,{\E}^{0,1}_2(\Rb^1\pi_{\d *}{\OO}_{\Xs}^*)).
\end{equation}
\indent
Consider the commutative diagram with exact rows
$$
\begin{CD}
0 @>>> 
\H^1(X,{\R}^0{\pi_{0*}}{\OO}^*_{X_0}) @>>> \Pic(X_0) @>>>
\H^0(X,\R^1\pi_{0*}{\OO}^*_{X_0})\\
@. @V{\alpha}VV @VVV @VVV\\
0 @>>>
\H^1(X,{\R}^0{\pi_{1*}}{\OO}^*_{X_1}) @>>> \Pic(X_1) @>>>
\H^0(X,\R^1\pi_{1*}{\OO}^*_{X_1})
\end{CD}
$$ 
induced by the Leray spectral sequences of $\pi_0$ and $\pi_1$. 
Evidently the quotient $\ker\{\Pic(X_0)\rightarrow \Pic(X_1)\}/
\ker \alpha$ injects into 
$\H^0(X,{\E}^{0,1}_2(\Rb^1\pi_{\d *}{\OO}^*_{\Xs}))$. 
\\
\indent
On each component we have an exact sequence analogous to 
(\ref{F11exact}), hence
$$\ker \alpha = \ker\{\H^1(X,{\R}^0{\pi_{0*}}{\cF}^{1,1}_{0}) 
\rightarrow \H^1(X,{\R}^0{\pi_{1*}}{\cF}^{1,1}_{1})\}.$$
Let ${\cal C}=\coker\{{\Rb}^0\pi_{\d *}{\cF}^{1,1}_{\d}\rightarrow 
{\R}^0\pi_{*0}{\cF}^{1,1}_{0}\}$ and let ${\cal Q}= \coker\{{\cal C} 
\rightarrow {\R}^0\pi_{*1}{\cF}^{1,1}_{1}\}$. The sheaves $\cal C$ 
and $\cal Q$ are subsheaves of ${\R}^0\pi_{*1}{\cF}^{1,1}_{1}$ and 
${\R}^0\pi_{*2}{\cF}^{1,1}_{2}$ respectively and have no global 
sections by Lemma \ref{F11prop}. It follows that
\begin{equation}\label{keralphaF11}
\ker \alpha = \H^1(X, {\Rb}^0\pi_{\d *}{\cF}^{1,1}_{\d}).
\end{equation}
Consider the commutative diagram with exact rows
$$
\minCDarrowwidth7pt
\begin{CD}
0 @>>> \H^1(X,{\Rb}^0\pi_{\d *}{\OO}^*_{\Xs}) @>>> \Pic(\Xs) 
@>>> \H^0(X,{\Rb}^1\pi_{\d *}{\OO}^*_{\Xs})\\
@. @VVV @V{\delta}VV @VV{\gamma}V\\
0 @>>> \H^1(X,{\Rb}^0\pi_{\d *}{\cF}^{1,1}_{\d}) @>>> 
\Hb^1(\Xs,{\cF}^{1,1}_{\d}) @>>> 
\H^0(X,{\Rb}^1\pi_{\d *}{\cF}^{1,1}_{\d})
\end{CD}
$$
Here $\gamma $ is injective by Lemma \ref{R1inj}. Using Lemma \ref{F11prop} 
we can identify the middle vertical with the map $\Pic(\Xs)\rightarrow 
\ker\{(\Pic(X_0)\rightarrow \Pic(X_1)\}$ induced by the componentwise 
spectral sequence. Thus $\ker \delta =\T(\Xs)$ and $\coker \delta=\F$ is a 
finite group. Since $X$ is seminormal, $\Rb^0\pi_{\d *}{\OO}^*_{\Xs}=
{\OO}^*_X$ and $\H^1(X,{\Rb}^0\pi_{\d *}{\OO}_{\Xs}^*)\cong \Pic(X)$. 
An obvious diagram chase, combined with (\ref{keralphaF11}), gives an exact 
sequence
\begin{equation}
\label{cokernelfinite}
0\rightarrow \Pic(X)/\T(\Xs) \rightarrow \ker \alpha \rightarrow \F,
\end{equation}
where $\F$ is a finite group. Since $X$ is seminormal 
${\pi}^*_{\d}: \Pic^0(X)\cong \Pic^0(\Xs)$ by Lemma \ref{exponentialseq}. 
The above diagram shows the restriction of this map to the torus $\T(\Xs)
\subseteq \Pic^0(\Xs)$ is an isomorphism. Hence we have inclusions
\begin{equation}
\label{PmodTinclusions}
\Pic(X)/\T(\Xs)\subseteq \Pic(\Xs)/\T(\Xs) \subseteq 
\ker\{\Pic(X_0)\rightarrow \Pic(X_1)\}, 
\end{equation}
where the second map is induced by the componentwise spectral sequence
computing $\Pic(\Xs)$. From (\ref{cokernelfinite}) and (\ref{PmodTinclusions}) 
we see that the induced map
\begin{equation}\label{quotinj}
\coker\{\Pic(X)/\T(\Xs)\rightarrow \Pic(\Xs)/\T(\Xs)\} \rightarrow
\H^0(X,{\E}^{0,1}_2(\Rb^1\pi_{\d *}{\OO}^*_{\Xs}))
\end{equation}
has finite kernel. The source of this map is isomorphic to 
$\NS(\Xs)/\NS(X)$ and torsion free by Lemma \ref{exponentialseq}. Thus 
(\ref{quotinj}) is injective and we have (\ref{NSquotientinj}). 

We claim the injective map (\ref{NSquotientinj}) fits into the commutative 
square
\begin{equation}\label{dia3}
\minCDarrowwidth7pt
\begin{CD}
\NS(\Xs)/\NS(X) @>{\text{into}}>>  \ \ \ \ 
\H^0(X, {\E}^{0,1}_2(\Rb^1\pi_{\d *}{\OO}^*_{\Xs}))\\
@V{\delta}VV @VVV\\
\Hb^1_{\zar}(\Xs, {\cH}^1_{\d})/\im \H^2_{\zarloc}(X,\Z)
@>>> \H^0_{\zar}(X,{\Rb}^1\pi_{\d *}{\cH}^1_{\d})
\end{CD}
\end{equation}
To explain the bottom row in the above diagram, we note that the spectral sequence 
of $X_{\an}\rightarrow X_{\zar}$ induces a (surjective) map
\begin{equation}\label{H2loc}
\H^2_{\zarloc}(X,\Z)\rightarrow \ker\{\H^1_{\zar}(X,{\cal H}^1_X)\rightarrow 
\H^3_{\zar}(X,{\cal H}^0_X)\}.
\end{equation}
This map fits into the commutative square of mixed Hodge structures
$$
\begin{CD}
\minCDarrowwidth6pt
\H^2_{\zarloc}(X,\Z) @>>> \H^1_{\zar}(X,{\cal H}^1_X)\\
@VVV @VVV\\
\Hb^1_{\zar}(\Xs, {\cH}^1_{\d}) @>=>> \Hb^1_{\zar}(\Xs, {\cH}^1_{\d})
\end{CD}
$$
where the kernel of the left vertical map is ${\W}_0$. Thus we have an 
injective map
$$\H^2_{\zarloc}(X,\Z)/{\W}_0 \rightarrow  
\H^1_{\zar}(X,{\cal H}^1_X)/\K,$$
where $\K=\ker\{\H^1_{\zar}(X,{\cal H}^1_X) \rightarrow 
\Hb^1_{\zar}(\Xs, {\cal H}^1_{\d})\}$. 
Now the bottom row of (\ref{dia3}) is obtained from the Leray spectral
sequence of $\pi_{\d}$ using that the sheaf map 
${\cal H}^1_X \rightarrow {\Rb}^0\pi_{\d *}{\cH}^1_{\d}$ 
induces an injection of $\H^1_{\zar}(X,{\cal H}^1_X)/\K$ into 
$\H^1_{\zar}(X,{\Rb}^0\pi_{\d *}{\cH}^1_{\d})$.  

Next we show the right vertical map in (\ref{dia3}) is injective. Since the quotient
$\NS(\Xs)/\NS(X)$ is torsion free, it suffices to show the following map is injective
\begin{equation}\label{gamma}
\gamma:
\H^0(X,{\E}^{0,1}_2(\Rb^1\pi_{\d *}{\OO}^*_{\Xs}))_{\Q}
\rightarrow 
\H^0(X, {\E}^{0,1}_2(\Rb^1\pi_{\d *}{\cH}^1_{\d}))_{\Q}.
\end{equation}
Let $\Pic(X_0)_{X}$ be the constant sheaf $\Pic(X_0)$ on $X$. We have 
an exact sequence
$$\Pic(X_0)_{X}\rightarrow \R^1\pi_{0*}\OO^*_{X_0}\rightarrow 
\R^1\pi_{0*}{\cH}^1_{X_0}
\rightarrow 0$$
and we define ${\cF}^0$ to be the image of $\Pic(X_0)_X$ in 
$\R^1\pi_{0*}\OO^*_{X_0}$. Starting with
the map $\Pic(X_1)_{X}\rightarrow \R^1\pi_{1*}\OO^*_{X_1}$, one defines 
analogously ${\cF}^1$. Then, by definition, 
$$\ker \gamma = \ker\{\H^0(X,{\cF}^0)_{\Q}\rightarrow 
\H^0(X, {\cF}^1)_{\Q}\}.$$ 
We show $\ker\{{\cF}^0\rightarrow {\cF}^1\}_{\Q}=0$. Let $x\in X$ be any 
point. For $i=0,1$ define $\F_{x}^i=(\pi_i^{-1}(x)_{\operatorname{red}})
_{\operatorname{sn}}\subseteq X_i$,~i.e.~${\F}^i_{x}$ is the 
seminormalization of the fiber of $x$ over $\pi_i$ (with its reduced 
structure). Set
$(X_i)_x=\Spec({\OO}_{X,x})\times_X X_i$. By definition, the stalk of 
${\cF}^i$ at $x$ is the image of $\Pic(X_i) \rightarrow \Pic((X_i)_x)$
for $i=0,1$. By Lemma \ref{piclemma} the map ${\cF}^i_x\rightarrow
\Pic({\F}^i_x)$ is injective, and clearly factors through 
$\Pic^0({\F}^i_x)$. We have so the following commutative square
$$
\minCDarrowwidth7pt
\begin{CD}
{\cF}^0_x @>{\text{into}}>> \Pic^0({\F}^0_{x})\\
@VVV @VVV \\
{\cF}^1_x @>>> \Pic^0({\F}^1_{x})
\end{CD}
$$
Now ${\F}^{\d}_x \rightarrow \{x\}$ is a proper hypercovering. 
Hence ${\E}_{\infty}^{1,0}=\ker\{\H^1({\F}^0_x,\Z)\rightarrow \H^1({\F}^1_x,\Z)\}$
is a quotient of $\Hb^1({\F}^{\d}_x,\Z)=\H^1(\{x\},\Z)=0$. This implies
the map $\Pic^0({\F}^0_{x})\rightarrow \Pic^0({\F}^1_{x})$ has finite kernel, 
and $\ker\{{\cF}^0\rightarrow {\cF}^1\}_{\Q}=0$ as claimed.  
\qed \end{proof}

\end{document}